# Shifted Brownian Fluctuation Game

SONG-KYOO (AMANG) KIM

## ABSTRACT

This article analyzes the behavior of a Brownian fluctuation process under a mixed strategic game setup. A variant of a compound Brownian motion has been newly proposed as called Sifted Brownian Fluctuation Process to predict turning points of a stochastic process. The Shifted Brownian Fluctuation Game has been constructed based this new process to find the optimal moment of actions. Analytically tractable results are obtained by using the fluctuation theory and the mixed strategy game theory.

**Keywords:** Brownian motion process; fluctuation theory; mixed game strategy; shifted Brownian fluctuation game;

## 1. INTRODUCTION

This research proposes an alternative variant of a one-dimensional Wiener process which can describes the random positions with containing ups and downs. The Sifted Brownian Fluctuation Process (SBFP) which is a compound Wiener process with state dependent conditions has been designed to predict turning points of a stochastic process. The SBFP which combined with the first exceed theory is able to find the first moment of a turning point either a concave or a convex shape. The first exceed theory is the compound process evolves until one of its marks hits (i.e. reaches or exceeds) its associated level for the first time and the process will evolve until one of the components hits its assigned level for the first time [30]. The first exceed theory delivers a closed joint functional to predict the moment of the first observed threshold which is crossing a turning point [31, 32] of the SBFP.

A mixed strategy is an assignment of a probability to each pure strategy [37]. When enlisting mixed strategy, it is often because the game doesn't allow for a rational description in specifying a pure strategy for the game. This allows for a player to randomly select a pure strategy [37]. Since probabilities are continuous, there are infinitely many mixed strategies available to a player. Since probabilities are being assigned to strategies for a specific player when discussing the payoffs of certain scenarios the payoff must be referred to as expected payoff. A mixed game strategy could be constructed on the top of a SBFP which represents random changes including economic changes, oil market changes and stock market changes. This two-players game is targeted to find a payoff matrix when the decision is made at one step prior to hit the first turning point. The Shifted Brownian Fluctuation Game (SBFG) is the two-person mixed strategy game with parameters in a payment matrix from a SBFP.

## 2. SHIFTED BROWNIAN FLUCTUATION GAME

### 2.1. Shifted Brownian Fluctuation Process

The Sifted Brownian Fluctuation Process (SBFP) is a compound Brownian motion process with state dependent conditions. Before formulating the SBFP, a simple Brownian motion process $Y(s)$ is denoted as follows:

$$Y(s) = \chi \sum_{k=1}^{\infty} X_k \cdot \mathbf{1}_{\left\{k \leq \left\lfloor \frac{s}{\Delta s} \right\rfloor\right\}} \tag{2.1}$$

where

$$X_k = \begin{cases} +1, & \text{$k$-th step is to length $\chi$ is to up,} \\ -1, & \text{$k$-th step is to down} \end{cases} \tag{2.2}$$

and $X_k$ are assumed independent with a fair probabilities:

$$P\{X_k = 1\} = P\{X_k = -1\} = \frac{1}{2}. \tag{2.3}$$

It is noted that $Y(s)$ is normally distributed with zero mean and variance $\sigma^2 t$ where $\chi^2 = \sigma^2 \cdot \Delta s$. From (2.1), the variant of the Brownian motion process $a(s)$ (aka. SBFP) as follows:

$$a(s) = w(\text{s}) + Y(s), \tag{2.4}$$

where $w$ is the mean changes between $[0, s)$ when a Brownian motion process is moving with time $s$. This process is basically same as a Brownian motion process except for "shifting" the mean in timely manner with the slope $w$ and $W(s) \sim \mathcal{N}(a^0 s, \sigma^2 s)$. All processes for the SBFP are defined on a probability space $(\Omega, \mathcal{F}, P)$ and $\mathcal{F}_A, \mathcal{F}_\tau \subseteq \mathcal{F}(\Omega)$ are $\sigma$-subalebras [30, 31, 37]. From (2.4), we have

$$W_n := a(\Delta_n) = w_n \Delta_n + Y(\Delta_n). \tag{2.5}$$

and suppose

$$\mathcal{A} := \sum_{n=0}^{\infty} W_n \varepsilon_{s_n}, s_0(=0) < s_1 < \cdots, \text{a.s.} \tag{2.6}$$

becomes a $\mathcal{F}_A$-measurable variant Brownian motion process which follows the normal distribution with the mean $\sum_{n=0}^{\infty} a_n^0 (s_n - s_{n-1})$ and the variance

$\sigma^2 \sum_{n=0}^{\infty}(s_n - s_{n-1})$. This process is observed at random moments in accordance with the point process:

$$\mathcal{T} = \sum_{n=0}^{\infty} \varepsilon_{\tau_n}, \tau_0 = 0. \tag{2.7}$$

From (2.1) and (2.4)-(2.6), the Shifted Brownian Fluctuation Process is defined as follows:

$$A_\tau = \sum_{k \geq 0} W_k \varepsilon_{\tau_k} = A_0 + \sum_{k=1}^{n} \{w_k \Delta_k + Y(\tau_n)\} \tag{2.8}$$

where $w_k, k = 1, 2, \ldots, w_0 = 0$ are the statement dependent constant values of SBFP with the notation

$$\Delta_k := \tau_k - \tau_{k-1}, k = 0, 1, \ldots, \tau_{-1} = 0, \tag{2.9}$$

and the following functionals can be evaluated

$$\gamma(v, \theta) = \mathbb{E}\left[e^{-vW_k - \theta \Delta_k}\right], \text{Re}(v) \geq 0, \text{Re}(\theta) \geq 0. \tag{2.10}$$

By using double expectation, the marginal Laplace-Stieltjes transform is applied as follows:

$$\gamma_k(\omega, \theta) = \delta\left(w_k \omega + \left(\frac{\sigma^2}{2}\right)\omega^2 + \theta\right) \tag{2.11}$$

where

$$\delta(\theta) = \mathbb{E}\left[e^{-\theta \Delta_1}\right]. \tag{2.13}$$

Analogously,

$$\gamma_0(\omega, \theta) = \mathbb{E}\left[\mathbb{E}\left[e^{-\omega A_0 - \theta \tau_0}\Big|\tau_0\right]\right] = \delta_0\left(A_0 \omega + \left(\frac{\sigma^2}{2}\right)\omega^2 + \theta\right) \tag{2.14}$$

where

$$\delta_0(\theta) = \mathbb{E}\left[e^{-\theta \tau_0}\right]. \tag{2.15}$$

In the Sifted Brownian Fluctuation Model (SBFG), the SBFP is ended when $A_\tau$ passes the first turning point with the correspond time $t^*$. With $S = [0, t^*)$, $t^* \in \mathbb{R}^+$, we are focused the time of turning points upon its escape from $S$. To formalize this model, the exit index is introduced as follows:

$$\nu := \inf\{k : \tau_k \geq t^*\}, \text{Re}(t^*) > 0, \tag{2.16}$$

and $\tau_\nu$ is the exit time or first passage time and $A_\nu$ is the position of the fluctuation at $\tau_\nu$. The actual moment of hit the turning point is $t^*$ and the first exceed value could be either a maximum (a concave shape process which is monotony decreased before $t^*$) or a minimum (a convex shape process which is monotony increased before $t^*$). The associated exit time from the confined SBFP and the formula (2.7) will be modified as

$$\widehat{A}_\tau = \sum_{k=0}^{\nu} W_k \varepsilon_{\tau_k}, \qquad (2.17)$$

which is the path of the SBFP from $\mathcal{F}(\Omega) \cap \{W_1 > w_1, W_2 > w_2, \ldots, W_{\nu-1} > w_{\nu-1} \cap W_\nu < w_\nu\}$ for a concave shape process (or $\mathcal{F}(\Omega) \cap \{W_1 < w_1, W_2 < w_2, \ldots, W_{\nu-1} < w_{\nu-1} \cap W_\nu > w_\nu\}$ for a convex shape). It gives an exact definition of the process observed until $\tau_\nu$. The functional

$$\Phi_\nu = \Phi_\nu(u, v, \vartheta, \theta) \qquad (2.18)$$

$$= \mathbb{E}\left[e^{-uA_{\nu-1}} e^{-vA_\nu} e^{-\vartheta\tau_{\nu-1}} e^{-\theta\tau_\nu} \mathbf{1}_{\{W_1 > w_1, W_2 > w_2, \ldots, W_{\nu-1} > w_{\nu-1} \cap W_\nu < w_\nu\}}\right]$$

The latter is of particular interests, we are interested in the observation moment of passing the turning point and one observation prior to this. The Theorem-SBFP establishes an explicit formula for $\Phi_\nu$ from (2.17)-(2.19). The Laplace-Carson transform is applied as follows:

$$\widehat{\mathcal{L}}_h(\bullet)(x) = x \int_{h=0}^{\infty} e^{-xh}(\bullet)dh, \ \mathrm{Re}(x) > 0, \qquad (2.21)$$

with the inverse

$$\widehat{\mathcal{L}}_x^{-1}(\bullet)(h) = \mathcal{L}^{-1}\left(\bullet\frac{1}{x}\right), \qquad (2.22)$$

where $\mathcal{L}^{-1}$ is the inverse of the bivariate Laplace transform [30, 31, 37].

**Theorem-SBFP.** The functional $\Phi_\nu$ of the SBFP on $\sigma$-algebra $\mathcal{F}(\Omega) \cap \{W_1 > w_1, W_2 > w_2, \ldots, W_{\nu-1} > w_{\nu-1} \cap W_\nu < w_\nu\}$ satisfies the following formula:

$$\Phi_\nu = \widehat{\mathcal{L}}_h^{-1}\left\{(\psi_0 - \psi_1) - \frac{\gamma_0 \cdot \phi}{2(\varphi^2 - 2\varphi)}(\Gamma_0 - \Gamma_1)\right\}(h^*). \qquad (2.23)$$

where
$$\psi_0 = \gamma_0(v, \theta), \qquad (2.36)$$
$$\psi_1 = \gamma_0(v, \theta + x), \qquad (2.37)$$
$$\gamma_0 = \gamma_0(u + v, \vartheta + \theta + x), \qquad (2.38)$$
$$\phi = \gamma_{\nu-1}(u + v, \vartheta + \theta + x), \qquad (2.39)$$
$$\varphi = \widehat{\gamma}(u + v, \vartheta + \theta + x), \qquad (2.40)$$
$$\Gamma_0 = \gamma_\nu(v, \theta), \qquad (2.41)$$
$$\Gamma_1 = \gamma_\nu(v, \theta + x). \qquad (2.42)$$

The moment of the first turning point $h^*$ is found as follows:

$$h^* = \left\{ h : \frac{d}{dh}\left( \frac{\partial}{\partial u}\sigma(h; u, 0, 0, 0)\Big|_{u=0} \right) = 0 \right\}, \tag{2.43}$$

where

$$\sigma(h; u, v, \theta, \vartheta) = \widehat{\mathcal{L}}_h^{-1}\left\{ (\psi_0 - \psi_1) - \frac{\gamma_0 \cdot \phi}{2(\varphi^2 - 2\varphi)}(\Gamma_0 - \Gamma_1) \right\}(h). \tag{2.44}$$

The functional $\Phi_\nu$ contains all decision making parameters regarding this standard stopping game. The information includes the first moments of a turning point ($\tau_\nu$), the moment of one step prior to passing the highest peak ($\tau_{\nu-1}$) and so on. The information from the closed functional are as follows:

$$\mathbb{E}[A_{\nu-1}] = \lim_{u \to 0}\left(-\tfrac{\partial}{\partial u}\right)\Phi_\nu(u, 0, 0, 0), \tag{2.45}$$
$$\mathbb{E}[\tau_{\nu-1}] = \lim_{\theta \to 0}\left(-\tfrac{\partial}{\partial \theta}\right)\Phi_\nu(0, 0, \theta, 0), \tag{2.46}$$
$$\mathbb{E}[A_\nu] = \lim_{v \to 0}\left(-\tfrac{\partial}{\partial v}\right)\Phi_\nu(0, v, 0, 0), \tag{2.47}$$
$$\mathbb{E}[\tau_\nu] = \lim_{\theta \to 0}\left(-\tfrac{\partial}{\partial \theta}\right)\Phi_\nu(0, 0, 0, \theta), \tag{2.48}$$

$$\mathbb{E}[\nu] = \left\lfloor \mathbb{E}\left[\tfrac{\tau_\nu}{\Delta_1}\right] \right\rfloor. \tag{2.49}$$

## 2.2. Shifted Brownian Fluctuation Game

A mixed game strategy could be constructed based on a SBFP which represents random changes including economic changes, oil market changes and stock market changes. This two-players game is targeted to find a payoff matrix when the decision is made at one step prior to hit the first peak at $h^*$ (see Fig. 2). The Shifted Brownian Fluctuation Game (SBFG) is the two-person mixed strategy game with parameters in a payment matrix from SBFP.

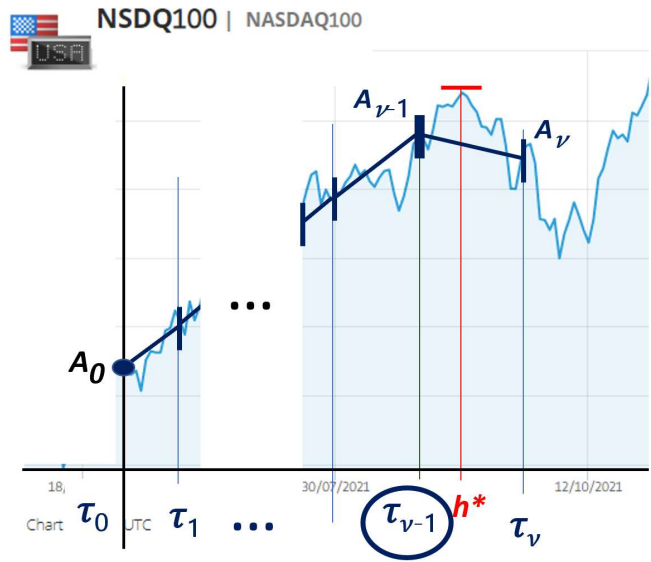

[Fig 2. SBFP on NASDAQ 100 stock chart for 6 months]

The players of this game is usually an uncontrollable subject (i.e., a nature, a market, an economy) verses a controller (i.e., a human, a company, a government). A controller (i.e., player 1) responses based on uncontrollable stochastic changes from an opponent player (i.e., player 2). In the SBFG, the decision is made at $\tau_{\nu-1}$ and the reward (payoff) of each player is determined after passing the peak (in a concave shaped SBFP). The normal of the game is:

$$\text{. Players:} \quad N = \{1, 2\}, \tag{2.50}$$
$$\text{. Strategy sets:} \quad s_1 = \{"Hold", "Action"\}, s_2 = \{"Up", "Down"\}. \tag{2.51}$$

## 3. Special Case: Memoryless Observation Process

Let us consider that the observation process has the memoryless property. This is very practical for actual implementation of the SBFG because this property implies that history of a SBFP is not considered. The moment of the decision making $h_{\nu-1}$ and the first exceed level index $\nu$ could be calculated from (2.49). Recalling from (2.22), the operator $\widehat{\mathcal{L}}_h$ is determined as follows:

$$G(x) = \widehat{\mathcal{L}}_h(f(h))(x), \tag{3.1}$$

and

$$f(h) = \widehat{\mathcal{L}}_x^{-1}(G(x)). \tag{3.2}$$

It is also noted that the formulas (2.11)-(2.15) could be rewritten as follows:

$$\delta(\theta) = \mathbb{E}\left[e^{-\theta \Delta_1}\right] = \left(1 + \widetilde{\delta} \cdot \theta\right)^{-1}, \tag{3.3}$$

$$\delta_0(\theta) = \mathbb{E}\left[e^{-\theta \tau_0}\right] = \left(1 + \widetilde{\delta}_0 \cdot \theta\right)^{-1}, \tag{3.4}$$

$$\widetilde{\delta} = \mathbb{E}[\Delta_1], \ \widetilde{\delta}_0 = \mathbb{E}[\tau_0] \tag{3.5}$$

$$\gamma_0(u, x) = \left(1 + \widetilde{\delta}_0 \cdot \left(\tfrac{\sigma^2}{2}u^2 + A_0 u + x\right)\right)^{-1} \tag{3.6}$$
$$\gamma_\nu(u, x) = \left(1 + \widetilde{\delta} \cdot \left(\tfrac{\sigma^2}{2}u^2 + w_\nu u + x\right)\right)^{-1} \tag{3.7}$$
$$\gamma_{\nu-1}(u, x) = \left(1 + \widetilde{\delta} \cdot \left(\tfrac{\sigma^2}{2}u^2 + w_{\nu-1} \cdot u + x\right)\right)^{-1} \tag{3.8}$$

$$\widehat{\gamma}(u, x) = \left(1 + \widetilde{\delta} \cdot \left(\tfrac{\sigma^2}{2}u^2 + \overline{w}u + \theta\right)\right)^{-1} \tag{3.9}$$
$$\overline{w} = \left(\tfrac{1}{k}\right)\sum_{j=1}^{k} w_k \tag{3.10}$$

Let us consider a monotonic increased SBFP (i.e., a concave shaped process) for this case. From (3.31), we can find

$$\mathbb{E}\left[e^{-uA_{\nu-1}}\right] = \Phi_\nu(u, 0, 0, 0) = \mathcal{L}_x^{-1}\Psi(x), \tag{3.11}$$

and

$$\Psi(x) = (1 - \gamma_0(0, x)) - \left(\frac{1}{2}\right)\left(\frac{\gamma_0(u,x)\gamma_{\nu-1}(u,x)\{1-\gamma_\nu(0,x)\}}{(\hat{\gamma}(u,x))^2 - 2\cdot\hat{\gamma}(u,x)}\right)$$

$$= G_0 + \frac{G_1}{2}\left(\frac{(1+D)^2}{(2D+1)(1+D_0)(1+D_2)}\right), \tag{3.12}$$

where

$$G_0 = \frac{\widetilde{\delta}_0 \cdot x}{\widetilde{\delta}_0 \cdot x + 1}, \tag{3.13}$$

$$G_1 = \frac{\widetilde{\delta} \cdot x}{\widetilde{\delta} \cdot x + 1}, \tag{3.14}$$

$$D = \widetilde{\delta} \cdot \left(\frac{\sigma^2}{2}u^2 + \overline{w}u + x\right), \tag{3.15}$$

$$D_0 = \widetilde{\delta}_0 \cdot \left(\frac{\sigma^2}{2}u^2 + A_0 u + x\right), \tag{3.16}$$

$$D_2 = \widetilde{\delta} \cdot \left(\frac{\sigma^2}{2}u^2 + w_{\nu-1} \cdot u + x\right). \tag{3.17}$$

From (2.43), (2.44) and (3.24), we can find the optimal moment of turning point $h^*$ by solving the following equation of $h$:

$$\frac{d\sigma^1(h)}{dh} = 0. \tag{3.27}$$

Therefore, we have

$$h^* = \left\{h \geq 0 : \left\{A \cdot e^{\left(\frac{1}{2\widetilde{\delta}}\right)h} = \frac{U - h}{h - V}\right\}, V \in \mathbb{R}^+\right\}, \tag{3.29}$$

where

$$U = \frac{(2+\widetilde{\delta})w_{\nu-1} - (3+\widetilde{\delta})\overline{w}}{\overline{w} - w_{\nu-1}}, \; U < 0, \tag{3.30}$$

$$V = \frac{(3+2\widetilde{\delta})\overline{w} - 2w_{\nu-1}}{\overline{w}}, \; V < 0, \tag{3.31}$$

$$A = \frac{\overline{w}}{2(w_{\nu-1} - \overline{w})}, \; A > 0. \tag{3.32}$$

From (3.30)-(3.31), we can find the condition of a first turning point as follows:

$$\left(\frac{3+\widetilde{\delta}}{2+\widetilde{\delta}}\right)\overline{w} \leq w_{\nu-1} \leq \left(\frac{3+2\widetilde{\delta}}{2}\right)\overline{w}. \tag{3.33}$$